\documentclass[letterpaper, 10 pt, conference]{ieeeconf}  

\IEEEoverridecommandlockouts                              
\usepackage[pdftex]{graphicx}
\usepackage{amsmath} 
\usepackage{amssymb}  
\usepackage[noadjust]{cite}
\usepackage{color}
\usepackage{enumerate}

\usepackage{balance} 

\usepackage{subfig}

\usepackage[colorlinks,hyperindex,pdfpagelabels=true,plainpages=false]{hyperref}
\usepackage{url}
\hypersetup{colorlinks,citecolor=black,filecolor=black,linkcolor=black,urlcolor=black}

\title{\LARGE \bf Practical Comparison of Optimization Algorithms for Learning-Based MPC with Linear Models}

\author{Anil Aswani$^\star$, Patrick Bouffard$^\star$, Xiaojing Zhang$^\star$, and Claire Tomlin
\thanks{This material is based upon work supported by NSF ActionWebs (	CNS-0931843), ONR MURI SMARTS (N00014-09-1-1051), and an NSERC fellowship (P. Bouffard). The views and conclusions contained in this document are those of the authors and should not be interpreted as representing the official policies of the NSF, the ONR, or the NSERC.} 
	\thanks{X. Zhang is with the Automatic Control Laboratory, Department of Electrical Engineering and Information Technology, Swiss Federal Institute of Technology Zurich (ETH Zurich), 8092 Zurich, Switzerland.}%
	\thanks{A. Aswani is with the Department of Industrial Engineering and Operations Research, University of California, Berkeley, CA 94720, USA.}%
		\thanks{P. Bouffard and C. Tomlin are with the Department of Electrical Engineering and Computer Sciences, University of California, Berkeley, CA 94720, USA. }%
\thanks{$^\star$ These authors contributed equally to this work.}
}

\begin{document}

\newtheorem{thm}{Theorem}
\newtheorem{cor}[thm]{Corollary}
\newtheorem{lem}[thm]{Lemma}

\maketitle
\thispagestyle{empty}
\pagestyle{empty}

\begin{abstract}
Learning-based control methods are an attractive approach for addressing performance and efficiency challenges in robotics and automation systems.  One such technique that has found application in these domains is learning-based model predictive control (LBMPC).  An important novelty of LBMPC lies in the fact that its robustness and stability properties are independent of the type of online learning used. This allows the use of advanced statistical or machine learning methods to provide the adaptation for the controller.  This paper is concerned with providing practical comparisons of different optimization algorithms for implementing the LBMPC method, for the special case where the dynamic model of the system is linear and the online learning provides linear updates to the dynamic model.  For comparison purposes, we have implemented a primal-dual infeasible start interior point method that exploits the sparsity structure of LBMPC.  Our open source implementation (called LBmpcIPM) is available through a BSD license and is provided freely to enable the rapid implementation of LBMPC on other platforms.  This solver is compared to the dense active set solvers LSSOL and qpOASES using a quadrotor helicopter platform.  Two scenarios are considered: The first is a simulation comparing hovering control for the quadrotor, and the second is on-board control experiments of dynamic quadrotor flight.  Though the LBmpcIPM method has better asymptotic computational complexity than LSSOL and qpOASES, we find that for certain integrated systems (like our quadrotor testbed) these methods can outperform LBmpcIPM.  This suggests that actual benchmarks should be used when choosing which algorithm is used to implement LBMPC on practical systems.
\end{abstract}

\section{Introduction}

There is growing interest in the development, application, and integration of learning-based control methods towards actuation and planning for robotic systems \cite{Tedrake2005,Lehnert2012}, unmanned autonomous vehicles \cite{abbeel2010,bouffard2011,aswani2012acc_quad,Keivan2012,Chowdhary2013}, and energy systems \cite{Aswani2011,aswani2012nmpc,Kolter2012}.  These techniques are able to achieve high performance and efficiency by learning better system models from measured data, which can then be used to more accurately control and plan the system.  Because the use of these methods requires online computation for both the model learning and the generation of a control sequence, optimized algorithms for both aspects of such learning-based control methods are important for practical implementations.

This paper focuses on algorithms for generating a control sequence when using the learning-based model predictive control (LBMPC) technique \cite{Aswani2011LBMPC}, restricted to the special case with linear dynamic models and linear learning to update the dynamic models.  (Note that the general LBMPC technique described in \cite{Aswani2011LBMPC} can handle nonlinear learning and nonlinear dynamic models.)  LBMPC uses statistical learning methods to improve the model of the system dynamics, while using robustness techniques from model predictive control (MPC) \cite{Garcia1988, Qin03asurvey, Mayne2000} to ensure that stability and system constraints are deterministically maintained \cite{Aswani2011LBMPC}.  It is similar to adaptive MPC \cite{Fukushima2007301,adetola2011}, which only work for specific types of model learning. 

The design of algorithms for LBMPC is particularly interesting because variants have been used for a variety of practical applications (e.g.\ \cite{Aswani2011,bouffard2011,aswani2012acc_quad,aswani2012nmpc,Lehnert2012,Keivan2012,Chowdhary2013}).  At its core, MPC is formulated as an optimization problem in which a cost function is minimized with respect to constraints on (a) the states and inputs of the system and (b) a model of the dynamics of the system. LBMPC is differentiated from MPC in that online learning is used to update the dynamics model, and careful structuring of the constraints of the corresponding optimization problem can be used to deterministically guarantee stability and robustness of the resulting controller \cite{Aswani2011LBMPC}. 

\subsection{Optimization Algorithms for MPC}

Because LBMPC is based on MPC, we first discuss different classes of optimization algorithms that can be used to compute a control sequence for MPC.  Depending on the type of system model (e.g.\ linear or nonlinear) and the form of the constraints, there are a variety of approaches to solving the optimization problem corresponding to an MPC controller.  The most efficient structure occurs when the system model is linear (or affine), the cost function is quadratic, and the constraints are linear; in this situation, the optimization problem is a convex quadratic program (QP), which can be solved relatively easily.  Such formulations are sometimes called QP-MPC problems.

Explicit MPC \cite{Borrelli2003,Bemporad2002ExplLGQ,Tondel01,Mariethoz2009} is one approach to solving a QP-MPC.  In this procedure, a lookup table that gives the optimal control as a function of the initial states is computed offline. It is well-suited for situations where the QP-MPC is time-invariant.  This time-invariant structure occurs, for instance, when the system model is not changing; however, such time-invariance is not the case for LBMPC where online learning updates the model.  Furthermore, the number of entries in the lookup table can grow exponentially with input and state dimensions. 

An alternative approach is to use a QP optimization solver (e.g.\ \cite{Gill1986LSSOL,SNOPT,Ferreau2008}) at each time step.  Specialized solvers that can exploit the sparsity of QP-MPC have computational complexity that scales linearly in the prediction horizon of the QP-MPC \cite{Wang2010FastMPC,Rao98IMP2MPC}, as opposed to non-sparse solvers that scale cubicly.  Other approaches include automatic code generation using problem-tailored solvers \cite{Mattingley2011CVXGEN,Domahidi2012fastMPC}, and combining explicit and online MPC \cite{Zeilinger2011}.

\subsection{Optimization Algorithms for LBMPC}

One advantage of LBMPC is that its formulation and safety properties are independent of the statistical method used; however, the numerical solver used to compute the LBMPC control is dependent on the form of the statistical method.  This is because the structure of the resulting optimization problem depends upon the type of statistical method used with LBMPC.  For instance, the optimization is non-convex when a nonparametric statistical method is used to do the learning \cite{Aswani2011LBMPC}.  On the other hand, if the system dynamics are linear, the state and input constraints are linear, the cost function is quadratic, and the statistical method provides linear model updates, then LBMPC can be described by a QP.  Such a QP-LBMPC formulation can be found in many engineering problems, including quadrotor flight control \cite{bouffard2011,aswani2012acc_quad} and energy-efficient building automation \cite{Aswani2011}.

It turns out that QP-LBMPC has sparsity structure similar to QP-MP. This mean that sparse solvers can be designed which computationally scale well. In particular, we have implemented a primal-dual infeasible start interior point method (PD IIPM) based on Mehrotra's predictor-corrector scheme \cite{Mehrotra1992PDIPM} in {C++}, and named this solver {LBmpcIPM}.  Our open source implementation (distributed via the BSD license) can be downloaded from \url{http://lbmpc.bitbucket.org/}.  For comparison, we consider two dense active set solvers (LSSOL v1.05-4 and qpOASES v3.0beta) \cite{Gill1986LSSOL,Ferreau2008}.  Because these solvers do not consider the sparsity of QP-LBMPC, their scaling properties are worse than LBmpcIPM.  That being said, for particular applications it may be the case that a dense solver such as LSSOL or qpOASES requires less computation than LBmpcIPM.  Our aim in this paper is to report simulation and experimental comparisons that provide insights into the practical implementation of algorithms for LBMPC.

We begin by formally defining the LBMPC technique in Section \ref{sec: lbmpc}. Next, we briefly summarize the key characteristics of the LBmpcIPM, LSSOL, and qpOASES algorithms that are being compared.  Section \ref{sec: simulations} presents numerical simulation results and real-time experiments.  A quadrotor helicopter testbed is used to provide a practical platform for comparison of different algorithms.  Two scenarios are considered: The first is a simulation comparing hovering control for the quadrotor, and the second is on-board control experiments of dynamic quadrotor flight.

\section{Learning-Based Model Predictive Control}
\label{sec: lbmpc}

For simplicity, we focus on QP-LBMPC, in which the cost is quadratic, constraints are polyhedral, and all dynamics are affine.  This results in a linear control law which is known to be computationally tractable and robust for practical applications.  The general form of LBMPC, which can handle nonlinear system dynamics and nonlinear learning, can be found in \cite{Aswani2011LBMPC}. 

\subsection{Three Models of System Dynamics}


To describe LBMPC, we must define three discrete-time models of the system. The first is the \textit{true  system} with dynamics
\begin{equation*}
	x_{m+1} = Ax_{m} + Bu_{m} + s + g(x_{m},u_{m}),
\end{equation*}
where $m $ denotes time, $x\in\mathbb{R}^n$ is the state, $u\in\mathbb{R}^m$ is the input, and $A\in\mathbb{R}^{n\times n},\ B\in\mathbb{R}^{n\times m},\ s\in\mathbb{R}^{n}$.  Here, $g:\mathbb{R}^{n}\times\mathbb{R}^{m}\to\mathbb{R}^n$ represents (possibly nonlinear) unmodeled dynamics.  The second model (\textit{nominal model}) is affine
\begin{equation}
\label{eq: affineSystem}
	\bar{x}_{m+1} = A\bar{x}_{m} + B\bar{u}_{m} + s,
\end{equation}
where $\bar{x}_m\in\mathbb{R}^n$, $\bar{u}_m\in\mathbb{R}^m$ are the state and control of the nominal model. The nominal model ensures two things: first, it is used to guarantee deterministic stability. Second, it is the model on which the learning is based, as discussed below.

The third model is the \textit{learned model}.  The learning in LBMPC occurs through a function $\mathcal{O}_m(\tilde{x}_{m},\tilde{u}_{m})$ that is known as an oracle.   The oracle provides updates to the nominal model in the following manner:
\begin{equation}
\label{eq: genOracleSystem}
	\tilde{x}_{m+1} = A\tilde{x}_{m} + B\tilde{u}_m + s + \mathcal{O}_m(\tilde{x}_{m},\tilde{u}_{m}),
\end{equation}
where $\tilde{x}\in\mathbb{R}^n$, $\tilde{u}\in\mathbb{R}^m$ are the state and input of the oracle system.  The oracle provides corrections to the nominal model by learning the unmodeled dynamics $g$ online.  One of the notable characteristics of the LBMPC formulation is that its deterministic stability and robustness properties are independent of both the statistical method used to estimate the model updates and the mathematical structure of the oracle.

However, the mathematical structure of the oracle does affect the structure of the optimization problem that must be solved to compute the control action of LBMPC.  In general, the optimization problem will be non-convex, which can be difficult to solve numerically in real time.  However, if the oracle has an affine form, then the optimization problem describing LBMPC is a QP.  More specifically, we consider QP-LBMPC where the oracle is given by
\begin{equation*}
\label{eq: affOracle}
	\mathcal{O}_m(\tilde{x}_{m},\tilde{u}_{m}) = L_m\tilde{x}_{m} + M_m\tilde{u}_{m} + t_m,
\end{equation*}
where $L_m\in\mathbb{R}^{n\times n},\ M_m\in\mathbb{R}^{n\times m},\ t_m\in\mathbb{R}^n$ are time-varying and constantly updated by an appropriate parametric statistical method \cite{aswani2012acc_quad}.

\subsection{QP-LBMPC}
We can state the QP-LBMPC control as the solution to the below optimization problem.  The interpretation of this optimization is given in the next subsection.
\begin{align}\label{eq: LBMPC}
	\min_{C,\theta} \quad &  (\tilde{x}_{m+N|m}- x^\star_{m+N|m})^T \tilde{Q}_f (\tilde{x}_{m+N|m}- x^\star_{m+N|m}) {\ }+ 	\nonumber \\	
			& \sum_{i=0}^{N-1} \Big{\lbrace}
						(\tilde{x}_{m+i|m} - x^\star_{m+i|m})^T \tilde{Q} (\tilde{x}_{m+i|m} - x^\star_{m+i|m}) {\ }+ \nonumber \\
			& (\check{u}_{m+i|m} - u^\star_{m+i|m})^T R (\check{u}_{m+i|m} - u^\star_{m+i|m})  \Big{\rbrace} \ \\	
	\text{s.t.} \quad & \tilde{x}_{m|m} = \hat{x}_m,\quad \bar{x}_{m|m} = \hat{x}_{m} \nonumber, \\
		&  \tilde{x}_{m+i|m}= \tilde{A}_m \tilde{x}_{m+i-1|m} + \tilde{B}_m \check{u}_{m+i-1|m} + \tilde{t}_m, \nonumber \\
	 &	\bar{x}_{m+i|m}= A \bar{x}_{m+i-1|m} + B \check{u}_{m+i-1|m} + s, \nonumber\\
	 & \check{u}_{m+i-1|m} = K\bar{x}_{m+i-1|m} + c_{m+i-1|m},  \nonumber \\
	 & F_{\bar{x},i} \bar{x}_{m+i|m} \leq f_{\bar{x},i},\qquad i = 1,\ldots,N  \nonumber \\
	 & F_{\check{u},i} \check{u}_{m+i|m} \leq f_{\check{u},i}, \qquad i=0,\ldots,N-1 \nonumber \\
	 & F_{\bar{x}\theta} \bar{x}_{m+j|m} + F_\theta \theta \leq f_{\bar{x}\theta} \nonumber
\end{align}
for a single fixed value of $j$, such that $j\in\{1,\ldots,N\}$. $N$ is the prediction horizon, $C$ the vector containing $\{c_{m+i|m}\}_{i=0}^{N-1}$, $\hat{x}_{m}\in\mathbb{R}^n$ the initial state, $\tilde{x}_{m+i|m},\bar{x}_{m+i|m}$ ($\check{u}_{m+i|m}$) the states (inputs) at time $m+i$ predicted at time $m$, and $\theta \in \mathbb{R}^m$ parameterizes the set of states that can be tracked by a steady-state control value. The dynamics matrices updated with the oracle are $\tilde{A}_m\triangleq A+L_m\in\mathbb{R}^{n\times n},\ \tilde{B}_m\triangleq B+M_m\in\mathbb{R}^{n\times m},\ \tilde{t}_m\triangleq s+t_m\in\mathbb{R}^n$, and $\{x^\star_{m+i|m}\}_{i=1}^N$ ($\{u^\star_{m+i|m}\}_{i=0}^{N-1}$) are the states (inputs) we want to track. We assume $\tilde{Q}=\tilde{Q}^T\in\mathbb{R}^{n\times n},\ \tilde{Q}_f=\tilde{Q}_f^T\in\mathbb{R}^{n\times n},\ R=R^T\in\mathbb{R}^{m\times m}$ and $\tilde{Q}=\tilde{Q}^T>0$, $R=R^T>0$, $\tilde{Q}_f=\tilde{Q}_f^T>0$.

\subsection{Interpretation of QP-LBMPC}

At an abstract level, the idea of LBMPC is to maintain two models, \eqref{eq: affineSystem} and \eqref{eq: genOracleSystem}, of the system within the optimization problem.  A cost function that depends on the states of the learned model $\tilde{x}_m$ and the control inputs $\check{u}_m$ is minimized.  The same control inputs are used to check that input and state constraints are satisfied when applied to the nominal model $\bar{x}_m$.  This is the reason that the optimization is formulated so that the constraints are applied to model \eqref{eq: affineSystem} but not \eqref{eq: genOracleSystem}, while both models have the same control input $\check{u}$. Furthermore, LBMPC robustifies the constraints to handle the mismatch between the models and the true system, represented by the polytopes defined by $(F_{\bar{x},i},f_{\bar{x},i})$ and $(F_{\check{u},i}, f_{\check{u},i})$ \cite{Aswani2011LBMPC}.  Finally, we note that LBMPC uses a tracking formulation of MPC that allows tracking to reference points parameterized by $\theta$.  This necessitates the use of a terminal constraint set $([F_{\bar{x}\theta}\ F_\theta], f_{\bar{x}\theta})$. Further details on its computation can be found in \cite{Aswani2011LBMPC}.  Under appropriate conditions \cite{Aswani2011LBMPC}, the LBMPC formulation ensures robust recursive feasibility, robust constraint satisfaction, and is robustly asymptotically stable (RAS).  

We require the constraint matrices $F_{\bar{x},i}\in\mathbb{R}^{l_{\bar{x}}\times n}$, $F_{\check{u},i}\in\mathbb{R}^{l_{\check{u}}\times m}$, $F_\theta\in\mathbb{R}^{l_\theta\times m}$ to be full column rank, where $l_{\bar{x}}, l_{\check{u}},l_\theta$ are the respective number of constraints. The full column rank condition holds, e.g.\ for box constraints. For simplicity of presentation, we keep the number of linear inequality constraints constant for all $i$. Our results, however, also hold even if the number of constraints vary over $i$.  The last condition in \eqref{eq: LBMPC} guarantees persistent feasibility by requiring any intermediate state at time $m+j$ to lie within an approximation of the maximal admissible disturbance invariant set \cite{Kolmanovsky1998}. Our formulation is meant to approximate the limit as $\theta$ approaches infinity. It is derived from the tracking formulation of \cite{Limon2010RobustTubeBasedMPC} as proposed in \cite{Aswani2011LBMPC}. The approximation was observed to deliver good performance and results both in simulation and experiments.  This variant of the tracking formulation for LBMPC maintains the robust constraint satisfaction and recursive feasibility properties from the general LBMPC formulation.  However, the RAS property for this variant of LBMPC still needs to be checked. The gain $K\in\mathbb{R}^{m\times n}$ is chosen such that $(A+BK)$ is Schur stable \cite{Aswani2011LBMPC}.

From the cost we infer that \eqref{eq: LBMPC} is a convex optimization problem. At each time step, a problem of form \eqref{eq: LBMPC} with current state $\hat{x}_{m}$ is solved. Clearly, the solution is a (non-trivial) function in $\hat{x}_{m}$. QP-LBMPC (and LBMPC and MPC in general) works as follows: Let $(\{\check{u}^\text{opt}_{m+i|m}\}_{i=0}^{N-1},\{\tilde{x}^\text{opt}_{m+i|m}\}_{i=1}^{N},\{\bar{x}^\text{opt}_{m+i|m}\}_{i=1}^{N},\theta^\text{opt})$ be the optimizer at time $m$. The QP-LBMPC policy takes $\check{u}_{m|m}=\check{u}^\text{opt}_{m|m}$ as the current control action and waits until the new measurement $\hat{x}_{m+1}$ becomes available at the next time step. A new problem is solved and the procedure is repeated for each time step. In this paper, we describe a method with computational complexity linear in the prediction horizon $N$. We therefore assume that the optimization problem is feasible with finite optimal value.

\section{Characteristics of Compared Algorithms}

The QP solvers that we use to compute the control sequence given by LBMPC \eqref{eq: LBMPC} have different features.  LBmpcIPM is a sparse primal-dual infeasible start interior point method (PD IIPM) \cite{Nocedal2000NumOpt} that exploits the particular sparsity structure of QP-LBMPC, and it is similar to earlier research in fast MPC \cite{Wang2010FastMPC,Rao98IMP2MPC,Domahidi2012fastMPC}.  The main difference to fast MPC is in the sparsity pattern of the matrices involved in the computational algorithm, and details on the implementation of LBmpcIPM can be found in \cite{Zhang2011}.  One salient feature of LBmpcIPM is that its computational complexity scales as $O(N(m+n)^3)$, which is linear in $N$.

In contrast, both LSSOL and qpOASES are dense active set solvers.  Strictly speaking, LSSOL is optimized for constrained least-squares (CLS) optimization problems, and QP-LBMPC can be reformulated as a CLS problem so as to work well with LSSOL \cite{aswani2012acc_quad}.  On the other hand, qpOASES incorporates homotopy methods to exploit the recursive nature of MPC (and LBMPC) optimization problems.  The computational complexity for both LSSOL and qpOASES when solving QP-LBMPC problems scales as $O(N^3(m+n)^3)$, cubic in $N$. This is worse scaling than LBmpcIPM, but for small values of $N$ it may be the case that LSSOL and qpOASES outperform LBmpcIPM.

\section{Experiments and Simulations}
\label{sec: simulations}

This section empirically compares {LBmpcIPM}  with two dense active set solvers (LSSOL v1.05-4 and qpOASES v3.0beta)  \cite{Gill1986LSSOL,Ferreau2008}.  We begin with simulations on a model of a quadrotor helicopter. Simulations show that the computation time of the sparse PD IIPM does indeed scale linearly in the prediction horizon $N$.  This is in contrast to the dense active set solvers whose computation times are shown to scale cubicly in $N$.  These simulations are followed by experiments on a quadrotor helicopter, in which three solvers are empirically compared. The solvers were run in real time using the computer on-board the quadrotor helicopter which is slow in comparison to a desktop computer. Before describing these results, we begin by summarizing the quadrotor helicopter model that is used.

It is important to note why we use LBMPC on the quadrotor helicopter.  The helicopter is subject to complex, nonlinear physics due to different aerodynamic effects, and many of these effects are state dependent.  The challenge from a practical standpoint is that modeling these effects is difficult. So we utilize learning to identify these effects in an automated manner as the helicopter is flying.  The learning in essence is identifying a time-varying linearization that describes the complex helicopter physics, and this enables improved flight performance in experiments \cite{bouffard2011,aswani2012acc_quad}.




For the results presented here, all solvers were compiled using the GCC v4.6.3 compiler, and the corresponding iteration counts for all solvers remained primarily between 2 to 4 iterations.  The termination criterion used for the {LBmpcIPM} solver was residues with norms below 1e-3.  Moreover, the MPC formulations that were used to control this model placed box constraints on all of the states and inputs.  

\subsection{Quadrotor Helicopter Model}

Here we summarize the basic features of the model. There are ten states $x \in \mathbb{R}^{10}$ that correspond to three positions, three velocities, two angles, and two angular velocities.  Strictly speaking, a quadrotor helicopter has three angles (yaw, pitch, and roll) that define its orientation; for simplicity we keep the yaw angle fixed, and this is the reason that we only have two angles in our state.  The helicopter has three inputs $u \in \mathbb{R}^3$ that correspond to commanded values of thrust, pitch, and roll.  Each discrete time step (and the corresponding values of $A$,$B$, and $s$) represents 25ms sampling period or 40Hz sampling frequency. Further details about the quadrotor model we use can be found in \cite{bouffard2011}.

\balance    


\subsection{Computational Scaling in Horizon Length}

We conducted a simulation in which the quadrotor was commanded to move from a height of 2m above the ground to a height of 1m.  The reason that LBMPC can be useful in similar scenarios is that complex aerodynamic behavior leads to a change in the thrust of the helicopter when it approaches the ground. LBMPC allows the designer to explicitly specify that the model can change in this manner and then leverage statistical tools to learn and compensate for the effect of this phenomenon \cite{bouffard2011, aswani2012acc_quad}.  We used nonzero values of $L_m,M_m$ to represent learning that has taken place and has identified changes in the helicopter physics.

The simulations were run on three different computers of varying computational power and architecture.  The horizon sizes $N$ were allowed to range from 5 steps to 240 steps, which represents 0.125s to 6s of horizon time because of the 25ms sampling period of the model. The first computer is a Lenovo T410 laptop, with a dual-core Intel Core i7 processor running at 2.67 GHz, with 4 MB of cache, 8 GB of RAM, and running a 64-bit Linux operating system. The second is a Dell Precision 390 desktop, with a dual-core Intel Core2 CPU running at 2.4 GHz, 4 MB of cache, 2 GB of RAM, and running 32-bit Linux. Finally the third computer is a small form factor computer-on-module (CoM) which runs onboard the quadrotors used in our laboratory testbed as reported in \cite{bouffard2011}. It is based on a single-core Intel Atom Z530 CPU running at 1.6 GHz, with 512 KB of cache, 1 GB of RAM, and runs 32-bit Linux. 

Table \ref{tab: performance} lists the simulation results for LBmpcIPM with varying horizon $N$ on all three computers. For a given CPU the solution time is clearly linear in the prediction horizon $N$. This is what is expected for a sparse interior point solver that exploits the special structure of an MPC problem, and stands in contrast to dense solvers that scale cubicly. The latter fact can be verified in Table \ref{tab: scaling}, which reports the solution times for LBmpcIPM, LSSOL, and qpOASES on the Core i7 platform for $N=5$ to $N=30$. As theoretically predicted, the solution times for the dense active set solvers scale cubicly in $N$.

When the helicopter is in flight, its processor must handle additional overhead due to processes like measurement communication and file storage.  Therefore, the solution times reported for the Atom CPU represent a lower bound on the solve times for when the helicopter is in flight.

\begin{table}[!t]
\caption{Average simulation solve times (ms) for \mbox{LBmpcIPM} on three different platforms.}
\begin{center}
\begin{tabular}{|c||c|c|c|} \hline
$N$ & i7 & Core2 & Atom \\ \hline\hline
$5$ & 0.933 & 2.0 & 9.8  \\
$10$ & 1.83 & 4.0 & 20.1  \\
$15$ & 2.667 & 6.0 & 30.3  \\
$30$ & 5.646 & 12.1 & 60.7  \\
$60$ & 10.751 & 24.5 & 121.9  \\
$120$ & 22.818 & 49.9 & 245.7  \\
$240$ & 50.145 & 105.4 & 491.7  \\
 \hline
\end{tabular}
\end{center}
\label{tab: performance}
\end{table}

\begin{table}[!t]
\caption{Average simulation solve times (ms) for \mbox{LBmpcIPM}, LSSOL, and qpOASES on Core i7 platform.}
\begin{center}
\begin{tabular}{|c||c|c|c|} \hline
$N$ & LSSOL & qpOASES & LBmpcIPM \\ \hline\hline
$5$ & 0.136 & 0.352 & 0.933  \\
$10$ & 0.433 & 1.241 & 1.83  \\
$15$ & 0.948 & 3.114 & 2.667  \\
$30$ & 4.89 & 18.163 & 5.646  \\
 \hline
\end{tabular}
\end{center}
\label{tab: scaling}
\end{table}

\subsection{Experimental Comparison}

LBMPC controllers implemented using {LBmpcIPM}, LSSOL, and qpOASES were compared on a quadrotor helicopter testbed.  This experiment is an interesting comparison of the different solvers because (i) the Intel Atom Z530 processor onboard the helicopter is slow in comparison to a desktop computer, (ii) the optimization problem to compute the control must be computed within 25ms to enable the real time control, and (iii) the quadrotor has constraints placed on its state and inputs that correspond to physical constraints such as not crashing into the ground.  We used a horizon of $N=5$ because this was the largest horizon in which all the solvers could reliably terminate their computations during the 25ms sampling period for computing the control value.  Note that this is in contrast to the horizon of $N = 15$ that was used with the LSSOL solver in past experiments applying LBMPC to the quadrotor \cite{bouffard2011,aswani2012acc_quad}. This means that the benefits of the linearly-scaling computational complexity are not apparent in these experiments, nor are the benefits of a longer MPC horizon. However, it is worth noting that the current on-board quadrotor computer (which dates from 2009) could be replaced with another of similar size and power requirements. It would be compatible with our quadrotor platform \cite{Bachrach2012}, yet with about an order of magnitude better performance \cite{Meier2012}, which we can predict would enable, for LBmpcIPM, horizons around $N=30$. An experiment was conducted in which the helicopter was commanded to, starting from a stable hover condition, go left 1m and then go right 1m. This was repeated 10 times in quick succession.   The learning used in \cite{bouffard2011,aswani2012acc_quad} was enabled for this experiment.  

\begin{figure}[htbp]
\centering
\includegraphics{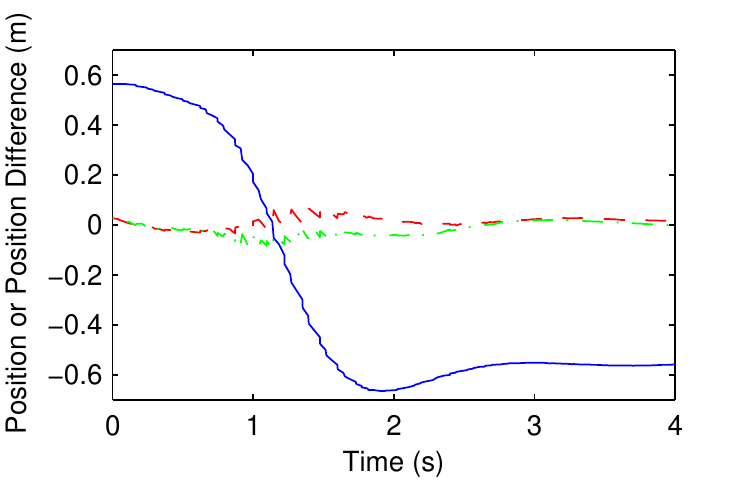}
\caption{The step response trajectory of the quadrotor helicopter flown using LBMPC solved with LSSOL is shown in solid blue,  the dashed red line indicates the difference between the trajectories of the helicopter when flown with the LSSOL versus the {LBmpcIPM} solver, and the dash-dotted green line denotes the difference between the trajectories of the helicopter when flown with the LSSOL versus the qpOASES solver.}
\label{fig:reprezent}
\end{figure}

A plot of a representative step input in which the quadrotor was commanded to go from left to right is shown in Fig.~ \ref{fig:reprezent}.  The position of the helicopter when using the LSSOL solver is shown in solid blue, the difference between the trajectories of the LSSOL and {LBmpcIPM} solvers is shown in dashed red, and the difference between the trajectories of the LSSOL and qpOASES solvers is shown in dash-dotted green.  We used the trajectory of the LSSOL solver as the reference trajectory, because this was the solver used for our previous experiments in \cite{bouffard2011,aswani2012acc_quad}.  As can be seen in the plots, the difference in trajectories is within 6.5cm for LBmpcIPM and 9.8cm for qpOASES. These differences are within a range that would be expected even between runs of the same trajectory using the same solver due to complex aerodynamic fluctuations that occur during a flight. They indicate that the different solvers are giving the same performance.

\begin{figure}
\centering
\subfloat[{LBmpcIPM} Solver]{\includegraphics{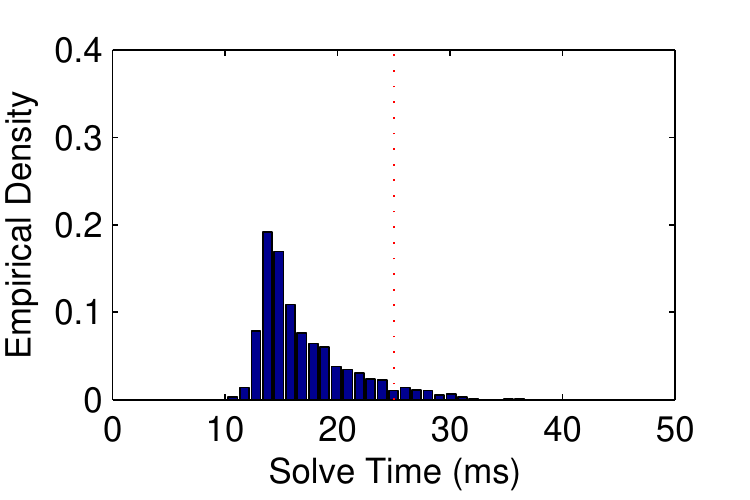} }\\
\subfloat[LSSOL Solver]{\includegraphics{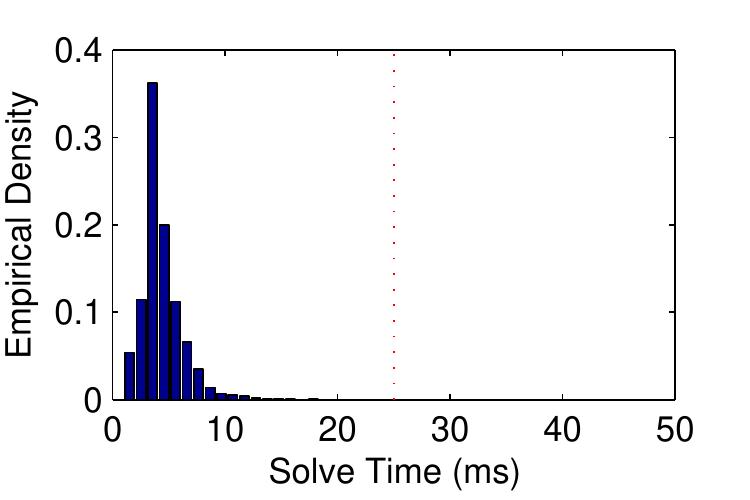} }\\
\subfloat[qpOASES Solver]{\includegraphics{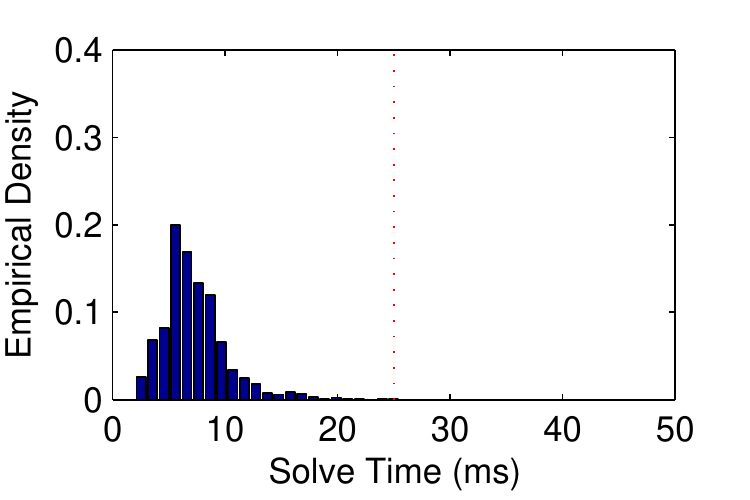} }
\caption[]{Empirical densities of solve times on quadrotor helicopter for different optimization algorithms are shown.  The vertical dashed red line at 25ms indicates the threshold beyond which greater solve times are too slow to be able to provide real time control.}
\label{fig:wassup}
\end{figure}

Histograms that show the empirical densities of the solve times for the three solvers can be seen in Fig.~\ref{fig:wassup}.  The two dense active set solvers have a lower variance of solve times, and this lower variance is important because it means that the solver is able to finish its computations under the time limit imposed by the 25ms sampling period.  In contrast, the \texttt{LBmpcIPM} method has a higher variance and so the solve times do exceed the 25ms limit a small percentage of the time.  We note that the advantage of our {LBmpcIPM} solver is that it has reduced computational effort, as compared to the dense active set solvers, when the horizon $N$ is large.  This benefit becomes readily apparent on faster computers that can handle longer horizons $N$, and this advantage of {LBmpcIPM} is not seen in our experiments where $N = 5$.

\section{Conclusions}
We have used simulations of quadrotor helicopter flight to confirm that the computation for a sparse solver like LBmpcIPM scales better than that of dense solvers like qpOASES and LSSOL.  However, real-time experiments of onboard implementations of these algorithms show that it may be the case that a dense solver like qpOASES or LSSOL can outperform the computational speed of LBmpcIPM.  This suggests that actual benchmarks should be used when choosing which algorithm is used to implement LBMPC on practical systems.  There is a last point that was hinted at in the paper but not explicitly discussed: The LBMPC problem has additional structure because of the similarity between the dynamics of the learned (with states $\tilde{x}$) and nominal model (with states $\overline{x}$), and this structure is not typical in linear MPC problems.  It may be possible to leverage this structure to provide improvements in the solve time.






\bibliographystyle{IEEEtran}
\bibliography{IEEEabrv,bibliography}

\end{document}